\documentclass[oneside,english]{amsart}
\usepackage[T1]{fontenc}
\usepackage[latin9]{inputenc}
\usepackage{babel}

\usepackage{verbatim}
\usepackage{amsthm}
\usepackage{graphicx}
\usepackage{amssymb}
\usepackage[numbers]{natbib}
\usepackage[unicode=true, pdfusetitle,
 bookmarks=true,bookmarksnumbered=false,bookmarksopen=false,
 breaklinks=false,pdfborder={0 0 1},backref=false,colorlinks=false]
 {hyperref}

\makeatletter
\numberwithin{equation}{section} 
\numberwithin{figure}{section} 
\theoremstyle{plain}
\newtheorem{thm}{Theorem}[section]
  \theoremstyle{remark}
  \newtheorem{rem}[thm]{Remark}
  \theoremstyle{plain}
  \newtheorem{lem}[thm]{Lemma}
  \theoremstyle{plain}
  \newtheorem{prop}[thm]{Proposition}
  \theoremstyle{plain}
  \newtheorem{cor}[thm]{Corollary}

\makeatother

\begin{document}

\title{IDLA on the Supercritical Percolation Cluster}
\begin{abstract}
We consider the internal diffusion limited aggregation (IDLA) process
on the infinite cluster in supercritical Bernoulli bond percolation
on $\mathbb{Z}^{d}$. It is shown that the process on the cluster
behaves like it does on the Euclidean lattice, in that the aggregate
covers all the vertices in a Euclidean ball around the origin, such
that the ratio of vertices in this ball to the total number of particles
sent out approaches one almost surely.
\end{abstract}

\author{Eric Shellef}

\address{The Weizmann Institute of Science\\
76100 Rehovot\\
Israel}

\email{shellef@gmail.com}



\date{Submitted to EJP on 18 August 2009, accepted on April 5, 2010}

\maketitle

\section{Introduction }

\subsection{Background and discussion}

Given a graph, IDLA defines a random aggregation process, starting
with a single vertex and growing by a vertex in each time step. To
begin the process, we specialize a vertex $v$ to be the initial aggregate
on the graph. In each time step, we send out one random walk from
$v$. Once this walk exits the aggregate, it stops, and the new vertex
is added to the aggregate. Let $\mathcal{I}(n)$ be the aggregate
in step $n$, thus $\mathcal{I}(1)=\left\{ v\right\} $ and $\left|\mathcal{I}(n)\right|=n$.

This process is a special case of a model Diaconis and Fulton introduced
in \citep{diaconis1991growth}. In the setting where the graph is
the $d$-dimensional lattice, Lawler, Bramson and Griffeath in \citep{lawler1992internal}
used Green function estimates to prove the process has a Euclidean
ball limiting shape. Let $B_{r}$ be all vertices in $\mathbb{Z}^{d}$
of Euclidean distance less than $r$ from the origin. The main theorem
in \citep{lawler1992internal} implies that for any $\epsilon>0$,
$\mathcal{I}(\left|B_{R}\right|)\supset B_{(1-\epsilon)R}$ for all
sufficiently large $R$ with probability one. Seeking to generalize
this result to other graphs, we note that convergence of a random
walk on the lattice to isotropic Brownian motion plays an important
role in convergence of IDLA on the lattice to an isotropic limiting
shape.

However, this property by itself is in general not enough for IDLA
to have such an inner bound, as the following example shows. Consider
the three dimensional Euclidean lattice, and choose a vertex $v$
at distance $R$ from the origin. Let $r=R^{0.9}$ and let $B_{r}(v)$
be the set of vertices of distance less than $r$ from $v$. Remove
all edges but one from the boundary of $B_{r}(v)$, and denote by
$v'$ the only vertex in $B_{r}(v)$ with a neighbor outside of $B_{r}(v)$.
Let us look at $\mathcal{I}(\left|B_{2R}\right|)$. A rough calculation
gives that the average number of visits to $v'$ is of order $R^{-1}\left|B_{2R}\right|$.
Since at least $\left|B_{r}(v)\right|=R^{2.7}$ visits to $v'$ are
needed to fill $B_{r}(v)$, we don't expect the ball to be full after
an order of $R^{3}$ particles have been sent out. Repeating this
edge removal procedure for balls of radius $R_{n}^{0.9}$ at distance
$R_{n}$ from the origin where $R_{n}=2^{n}R$, will ensure that there
is never a Euclidean inner bound. However, a random walk from the
origin on this graph will converge to Brownian motion because our
disruptions are sublinear. We will not give full proofs of these facts,
but hope they convince the reader that to get an inner bound some
kind of local regularity property is needed.

The main theorem of this paper states that IDLA on the supercritical
cluster in the lattice has a Euclidean ball as an inner bound. The
two main tools that are used to show this are a quenched invariance
principle \citep{berger2007quenched} which gives us convergence in
distribution to Brownian Motion from a fixed point, and a Harnack
inequality from \citep{barlow2004random}, which give us oscillation
bounds on harmonic functions in all small balls close to the origin.
The latter allows us to establish the local regularity missing from
the above example.%
\begin{figure}
\centering{}\includegraphics[scale=0.5]{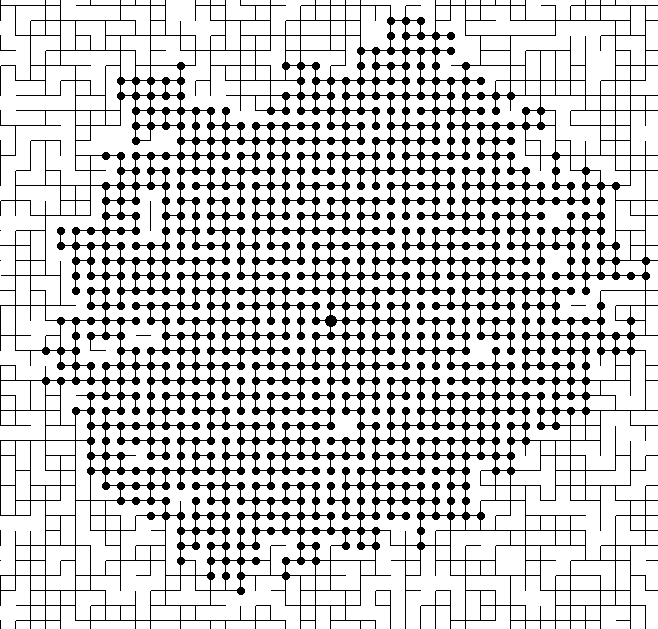}\caption{{\small IDLA of 1000 particles from marked vertex on the percolation
cluster with $p=0.7$.}}

\end{figure}

\subsection{\label{sub:Assumptions-and-statement}Assumptions and statement}

Consider supercritical bond percolation in $\mathbb{Z}^{d}$ with
the origin conditioned to be in the infinite cluster. Let the graph
$\Gamma(V,E)$ be the natural embedding in $\mathbb{R}^{d}$ of the
infinite cluster, i.e. $V\subset\mathbb{R}^{d}$. Fixing this embedding
for $\Gamma$, we get two separate (but comparable, see \citep{antal1996chemical})
distances. We denote by $\left|x-y\right|$ or $d(x,y)$ Euclidean
distance between points $x,y$, and by $d_{\Gamma}(x,y)$ the graph
distance between them. If one of the points is not in $V$, $d_{\Gamma}(x,y)=\infty$.
Let $B_{r}(x)=\{v\in V\ |\ d(x,v)<r\}$ be the vertices contained
in a Euclidean ball of radius $r$ and center $x$. We abbreviate
$B_{r}(\mathbf{0})$ as $B_{r}$. To differentiate between such a
set of vertices and a full ball in $\mathbb{R}^{d}$, we denote by
bold lowercase the following: $\mathbf{b}_{r}(x)=\{y\in\mathbb{R}^{d}\ |\ d(x,y)<r\}$.
Let $\mathbf{d}_{r}(x)$ be a box of side $r$ with center $x$, and
as above, let $D_{r}(x)$ be the vertices in this box. 

By existing work, which we will reference later, we know that with
probability one, the graph of the supercritical cluster, $\Gamma$,
in its natural $\mathbb{R}^{d}$ embedding, satisfies the following
assumptions:
\begin{enumerate}
\item \label{enu:conv_to_BM} A random walk on $\Gamma$ converges weakly
in distribution to a Brownian Motion as defined in \ref{sub:def-weak-conv-bm}.
\item \label{enu:Convergence-to-positive}Convergence to vertex density
$a_{G}>0$ as defined in \ref{sub:Convergence-to-positive}.
\item \label{enu:exit_time_uni_bnd}A uniform upper bound on exit time from
a ball as defined in \ref{sub:def-uni-bnd-exit-time}.
\item \label{enu:HK} A Harnack inequality as defined in \ref{sub:def-HK}.
\end{enumerate}
We denote by $X_{t}$ a discrete {}``blind'' random walk on $\Gamma$
defined as follows. For $x\in\Gamma$ and $t\in\mathbb{N}\cup\left\{ 0\right\} $,
$P\left(X_{t+1}=y|X_{t}=x\right)=1/2d$ if $y$ is a neighbor of $x$
in $\Gamma$, and $P\left(X_{t+1}=x|X_{t}=x\right)=1-\deg x/2d$.
We prefer this walk since its Green functions are symmetric, a fact
which will be useful later on. For non-integer $t$ we set $X_{t}=X_{\lfloor t\rfloor}$.

Note that assumption \ref{enu:HK} does not imply \ref{enu:exit_time_uni_bnd},
see \citep{delmotte2002graphs} for a counterexample.

Using only above assumptions and that $V\subset\mathbb{Z}^{d}$ (which
serves mostly to simplify notation and could be replaced by weaker
conditions), we show the IDLA process starting at $\mathbf{0}$ will
have a Euclidean ball inner bound as stated in the following theorem:

Let $\mathcal{I}(n)$ denote the random IDLA aggregate of $n$ particles
starting at $\mathbf{0}$.
\begin{thm}
\label{thm:IDLA-lower-bound} Almost surely, for any $\epsilon>0$,
we have that for all large enough $R$, $B_{(1-\epsilon)R}\subset\mathcal{I}(\left|B_{R}\right|)$ 
\end{thm}
We fix $\epsilon$ with which we prove the above for the rest of the
paper.

\subsection{Outline}

In the remainder of this section we state our assumptions precisely,
and explain why these assumptions are valid almost surely for the
infinite cluster in supercritical percolation.

Then, to prove IDLA inner bound, we need to show that for each vertex
$v\in V$ the Green function from \textbf{$\mathbf{0}$ }and expected
exit time from a ball around \textbf{0 }of a random walk starting
at $v$, behave similarly to those functions of a Brownian motion.
The exact statement needed appears in Lemma \ref{lem:pointwise_domination}.
The invariance principle gives us integral convergence of these functions
to the right value, but to improve them to pointwise statements, we
must show that they are locally regular. We use the Harnack inequality
from \citep{barlow2004random} to prove they are Hölder continuous.
This is a similar scheme to improving a CLT to an LCLT, as was done
in \citep{barlow2009parabolic} using a different method than ours.

In section \ref{sec:ptwise_bnds}, we start by proving a lemma comparing
expected exit time from a set to the Green function of a point in
the set. Then, we show how our Harnack assumption leads to an oscillation
inequality which we use to show regularity of the Green function and
expected exit time of the walk in a ball when we are far from the
boundary and the center. Next, in section \ref{sec:Domination-of-Green},
we use our assumption of an invariance principle to show that in small
balls the integral of these functions approaches a tractable limit
that can be calculated by knowledge of Brownian motion. Finally, in
section \ref{sec:Lower-Bound}, we utilize these estimates to prove
theorem \ref{thm:IDLA-lower-bound}.
\begin{rem}
Since the paper makes use of results based on ergodic theory, there
is no rate estimate for the convergence in theorem \ref{thm:IDLA-lower-bound}.
Also, on the Euclidean lattice there is an almost sure outer bound
with sublinear fluctuations from the sphere. Another interesting question
that is not treated here, is whether a similar outer bound holds in
our setting.
\end{rem}

\subsection{\label{sub:def-weak-conv-bm}Weak convergence to Brownian Motion}

In this and the next three subsections, we give a precise formulation
of each of our assumptions from above, and argue that they hold a.s.
on $\Gamma$, the supercritical cluster.

Let $\mathcal{C}_{T}^{d}=C([0,T]\to\mathbb{R}^{d})$, i.e. the continuous
functions from the closed interval $[0,T]$ to $d$-dimensional Euclidean
space. For $R\in\mathbb{N}$, let \[
w_{R}(t)=\frac{1}{R}\left(X_{\lfloor tR^{2}\rfloor}+(tR^{2}-\lfloor tR^{2}\rfloor)(X_{\lfloor tR^{2}\rfloor+1}-X_{\lfloor tR^{2}\rfloor})\right).\]
Thus $w_{R}(t)$ is a scaled linear interpolation of $X_{t}$ (defined
in \ref{sub:Assumptions-and-statement}) with its restriction to $[0,T]$
an element of $\mathcal{C}_{T}^{d}$.

We say that assumption \ref{enu:conv_to_BM} holds if for any $T>0,$
the law of $w_{R}(t)$ on $\mathcal{C}_{T}^{d}$ converges weakly
in the supremum topology to the law of a (not necessarily standard)
Brownian Motion $(B_{t}:0\le t\le T)$ . 
\begin{lem}
Assumption \ref{enu:conv_to_BM} holds for $\Gamma$ with probability
one.\end{lem}
\begin{proof}
This is Theorem 1.1 in \citep{berger2007quenched}. Another paper
with a similar invariance principle is \citep{mathieu2005quenched}. 
\end{proof}
Let $\mathcal{B}(t)$ denote the Brownian motion weak limit of $X_{t}$.

\subsection{\label{sub:Convergence-to-positive}Convergence to positive density}

Assumption \ref{enu:Convergence-to-positive} holds if there exists
a positive $a_{G}$ such that for any $\delta,\gamma>0$ and all sufficiently
large $R$ we have that for all $x\in\mathbf{b}_{R}$ and all $\delta R\le r\le R$
\[
\ \ \biggl|\frac{|D_{r}(x)|}{r^{d}}-a_{G}\biggr|,\ \biggl|\frac{|B_{r}(x)|}{\left|\mathbf{b}_{1}\right|r^{d}}-a_{G}\biggr|<\gamma.\]
That is, balls and boxes of size $\delta R$ to $R$, have vertex
density in $\left(a_{G}-\gamma,a_{G}+\gamma\right)$ for all large
enough $R$.
\begin{lem}
Assumption \ref{enu:Convergence-to-positive} holds for $\Gamma$
with probability one.\end{lem}
\begin{proof}
Let $\theta(p)$ be the probability for $\mathbf{0}$ to be in the
infinite cluster. $\theta(p)$ is positive in the supercritical regime.
From Theorem 3 in \citep{durrett1988large} for $d=2$ and from (2)
on p. 15 of \citep{gandolfi1989clustering} for $d>2$, we know that
for any $\rho>0$ there is a positive $c=c(\rho)$ such that \[
P\left(n^{-d}\left|D_{n}(x)\right|<(1-\rho)\theta(p)\right)<\exp(-cn).\]
where $x\in\mathbb{R}^{d}$ and $n\in\mathbb{N}$. Recall $D_{n}(x)$
are the vertices in a box of side-length $n$ and center $x$. Theorem
2 in \citep{durrett1988large} proves the easier density upper bound
for $d=2$, which under trivial modifications gives that for all $d$
\[
P\left(n^{-d}\left|D_{n}(x)\right|>(1+\rho)\theta(p)\right)<\exp(-cn).\]
The above, together with Borel Cantelli gives that if we choose $\gamma>0$,
$x\in\mathbf{b}_{1}$, then for all large enough $R$\[
\left|\left|D_{\rho R}(Rx)\right|-\theta(p)(\rho R)^{d}\right|<\frac{\gamma}{2}R^{d}.\]
So we have the result for small boxs of size $\rho R$. To expand
it to larger boxes, let $D$ be any box of diameter between $\delta R$
to $R$. We partition the box $D_{R}$ into $\rho R$-sized boxes
and choose $\rho=\rho(\delta)$ small enough so that the number of
boxes that intersect both $D$ and $D_{R}\backslash D$ is negligible
compared to the number that intersect $D$. Proving for balls is similar.
\end{proof}

\subsection{\label{sub:def-uni-bnd-exit-time}Uniform Bound on Exit time from
ball}

Denote by $\tau_{r}(x)$ the first time the walk leaves $B_{r}(x)$,
and write $\tau_{r}$ for $\tau_{r}(\mathbf{0})$.

Assumption \ref{enu:exit_time_uni_bnd} holds if there is a $c_{E}$
such that for any $\delta>0$ and all large enough $R$, if $x\in B_{R}$
and $r>\delta R$ then \begin{equation}
E_{x}\left[\tau_{r}(x)\right]<c_{E}r^{2}.\label{eq:E_x-pnt-up-bnd}\end{equation}
The assumption holds for $\Gamma$ with probability one as a direct
consequence of the following lemma, proved below.
\begin{lem}
With probability one, there is an $L$ such that for all $\delta>0$
and all large enough $R$, if $x\in B_{R}$ and $r>\delta R$ then\begin{equation}
P_{x}(\tau_{r}(x)>Lr^{2})<\exp(-cL).\label{eq:tau_R-uni-bdd-tail}\end{equation}
\end{lem}
\begin{proof}
Since we prove for all $R$ outside a bounded interval, it suffices
to prove the above for $r=\delta R$ with fixed $\delta>0$. To prove
this holds a.s. for supercritical percolation clusters, we use a heat
kernel upper bound given in Theorem 1 of Barlow's paper \citep{barlow2004random}.
The proof in \citep{barlow2004random} is for continuous walks with
a mean time of one between jumps. However, it can be transferred to
our discrete walk using Theorem 2.1 of \citep{berger2007quenched}.
We state an implication of what was proved.

With probability one there exists a function from $V$ to $\mathbb{N}$,
$\left\{ \mathcal{T}_{x}<\infty\right\} _{x\in V}$ where $\mathcal{T}_{x}$
is sublinear in the sense that $\frac{\mathcal{T}_{x}}{\|x\|}\to0$,
and there exist positive constants $c_{1},c_{2}$ such that for any
$x\in V$, the following holds. For all $y\in B_{t}(x),\ t\ge\mathcal{T}_{x}$:
\begin{equation}
P_{x}(X_{t}=y)<c_{2}t^{-d/2}.\label{eq:HK_up_bound}\end{equation}
Barlow's result actually states that almost surely $\mathcal{T}_{x}$
grows slower than a logarithmic function of $\|x\|$ and gives stronger
Gaussian bounds for the heat kernel from below and above. 

Using the heat kernel upper bound \eqref{eq:HK_up_bound}, the convergence
to zero of $\frac{\mathcal{T}_{x}}{\|x\|}$, and the upper bound on
vertex density resulting from $\Gamma\subset\mathbb{Z}^{d}$, we have
that for any $K$ and all large enough $R$, if $x\in B_{R}$ and
$r=\delta R$, then \begin{eqnarray*}
P_{x}(X_{Kr^{2}}\notin B_{r}(x)) & = & 1-\sum_{y\in B_{r}(x)}P_{x}(X_{Kr^{2}}=y)\\
 & \ge & 1-\sum_{y\in B_{r}(x)}c_{2}K^{-d/2}r^{-d}\\
 & \ge & 1-C(d)K^{-d/2}.\end{eqnarray*}
Thus for some $K'$, for all large enough $R$, and all $x\in B_{R}$,
$P_{x}(X_{K'r^{2}}\in B_{r}(x))<\frac{1}{2}$. Next, for any positive
$L$, we use the Markov property to upper bound $P_{x}(\tau_{r}(x)>LK'r^{2})$
by%
{}\[
P_{x}(X_{K'r^{2}}\in B_{r}(x))\prod_{i=2}^{\lfloor L\rfloor}P(X_{iK'r^{2}}\in B_{r}(x)\ |\ X_{(i-1)K'r^{2}}\in B_{r}(x))\le2^{-\lfloor L\rfloor},\]
 which proves the claim. 
\end{proof}
We will later use that assumption \ref{enu:exit_time_uni_bnd} implies
$L_{1}$ convergence of $R^{-2}\tau_{R}\wedge T$ to $R^{-2}\tau_{R}$
with $T$, uniformly in $R$ for all $x\in B_{R}$.
\begin{lem}
For all $\beta>0$ there is a $T(\beta)$ such that for all large
enough $R$, $\forall x\in B_{R}$\begin{equation}
\left|E_{x}\left[\tau_{R}\right]-E_{x}\left[\tau_{R}\wedge TR^{2}\right]\right|<\beta R^{2}.\label{eq:E_x-uni-conv-with-T}\end{equation}

\end{lem}
We rewrite the left hand side of \eqref{eq:E_x-uni-conv-with-T} and
use the Markov property with the exit time assumption. \begin{eqnarray*}
E_{x}\left[\left(\tau_{R}-TR^{2}\right)\mathbf{1}_{\left\{ \tau_{R}>TR^{2}\right\} }\right] & = & \sum_{y\in B_{R}}\biggl[E_{x}\left[\tau_{R}-TR^{2}\ |\ \tau_{R}>TR^{2},X(TR^{2})=y\right]\cdot\\
 &  & P(\tau_{R}>TR^{2},X(TR^{2})=y)\biggr]\\
 & = & \sum_{y\in B_{R}}E_{y}\left[\tau_{R}\right]P(\tau_{R}>TR^{2},X(TR^{2})=y)\\
 & \le & c_{E}R^{2}\sum_{y\in B_{R}}P(\tau_{R}>TR^{2},X(TR^{2})=y)\\
 & = & c_{E}R^{2}P(\tau_{R}>TR^{2}).\end{eqnarray*}
Since $E_{x}\left[\tau_{R}/R^{2}\right]$ is bounded by $c_{E}$ for
all $R$, we have by the Markov inequality that $P(\tau_{R}/R^{2}>T)\le c_{E}/T$
which converges to zero as $T$ goes to $\infty$ uniformly in $R$
for all $x\in B_{R}$.

\subsection{\label{sub:def-HK}Harnack inequality}

A function $h:\Gamma\to\mathbb{R}$ is harmonic on $x\in\Gamma$ if
$h(x)=\left|N(x)\right|^{-1}\sum_{y\in N(x)}h(y)$ where $N(x)=\left\{ y\in\Gamma:d_{\Gamma}(x,y)=1\right\} $
are the neighbors of $x$. We say $h$ is harmonic on $Z\subset\Gamma$
if it harmonic on every vertex $v\in Z$.

Assumption \ref{enu:HK} holds if there is a $1<c_{H}<\infty$ such
that for all $\delta>0$ we have that for all large enough $R$, if
$x\in B_{R}$, $r>\delta R$ and $h$ is a function that is non-negative
and harmonic on $B_{2r}(x)$ then \begin{equation}
\sup_{B_{r}(x)}h\le c_{H}\inf_{B_{r}(x)}h.\label{eq:HK}\end{equation}
In order to simplify the paper, the inequality above is formulated
for Euclidean distance rather than graph distance. We show that since
the distances are almost surely comparable on $\Gamma$, this is the
same. We start by proving the assumption holds for graph distance,
using Lemma 2.19 and Theorem 5.11 of \citep{barlow2004random}. We
write $B_{r}^{\Gamma}=\left\{ x\in\Gamma:d_{\Gamma}(\mathbf{0},x)\le r)\right\} $.
The lemma, using appropriate parameters and Borel Cantelli, tells
us that for all large enough $R$, $B_{R\log R}^{\Gamma}$ are \emph{very
good} balls. Specifically, all balls contained in $B_{R\log R}^{\Gamma}$
of graph radius larger than $R^{1/4(d+2)}$ have a positive volume
density and satisfy a weak Poincaré inequality as explained in (1.15)
and (1.16) of the same. Theorem 5.11 then tells us that for very good
ball of graph radius $R\log R$, all $x\in B_{\left(R\log R\right)/3}^{\Gamma}$
satisfy that for any function $h$ non-negative harmonic on $B_{R}^{\Gamma}(x)$
\begin{equation}
\sup_{B_{R/2}^{\Gamma}(x)}h\le\hat{c}_{H}\inf_{B_{R/2}^{\Gamma}(x)}h\label{eq:HK_graph}\end{equation}
where $\hat{c}_{H}(d,p)>0$ is a constant dependent only on dimension
and percolation probability. Since all but a finite number of balls
of graph radius $R\log R$ are very good and $R$ is $o(R\log R)$
we have assumption \ref{enu:HK} for graph distance.

Next, we transfer this to the Euclidean balls formulation of \eqref{eq:HK}.
By Theorem 1.1 of \citep{antal1996chemical}, we have that for some
$k>0$, $\rho(d,p)>0$ and $M<\infty$,\[
\mathbb{P}\left[d_{\Gamma}(x,y)>\rho m\Bigl|x,y\in\Gamma,\left|x-y\right|=m>M\right]<\exp(-km).\]
Let $A(x,y)$ be the event $\left\{ d_{\Gamma}(x,y)>\rho m\right\} $.
Union bounding the probability for $A(x,y)$ over every pair of points
in $B_{R}$ of Euclidean distance greater than $c\log R$ for some
large $c(k)$, we upper bound the probability of any such event occuring
by\[
CR^{d}\sum\limits _{m=c\log R}^{\infty}m^{d-1}\exp(-km)<C'R^{d}\left(\log R\right)^{d-1}R^{-2d},\]
which is summable for $d>1$. Hence by Borel Cantelli, almost surely
for all large enough $R$, we have that for any $x,y\in B_{R}$ with
$\left|x-y\right|>c\log R$, $d_{\Gamma}(x,y)<\rho\left|x-y\right|$.
Since $\log R$ is $o(R)$, this implies that for any $\delta>0$
and all large enough $R$, for any $x\in B_{R}$, \begin{equation}
B_{\delta\rho^{-1}R}(x)\subset B_{\delta R}^{\Gamma}(x).\label{eq:dist_comp}\end{equation}
Note that $B_{r}^{\Gamma}(x)\subset B_{r}(x)$ always because $\Gamma$
is embedded in $\mathbb{Z}^{d}$. Hence, given a function that is
non-negative harmonic on $B_{2r}(x)$ it is also such on $B_{2r}^{\Gamma}(x)$.
For large enough $R$, we use \eqref{eq:HK_graph} and \eqref{eq:dist_comp}
to get \eqref{eq:HK} for $B_{\rho^{-1}r}(x)$ with constant $\hat{c}_{H}$.
A routine chaining argument (see e.g. (3.5) in \citep{delmotte2002graphs})
transfers this to $B_{r}(x)$ as required with new constant $c_{H}$.

\section{\label{sec:ptwise_bnds}Pointwise bounds on Green function and expected
exit time}

In this section, we show the assumptions of a uniform bound on exit
time from a ball along with the Harnack inequality give us pointwise
bounds of the Green function and expected exit time of a random walk. 

As a convention, we use plain $c$ and $C$ to denote positive constants
that do not retain their values from one expression to another, as
opposed to subscripted constants $c_{i},$ that do. In general, these
constants are graph dependent, and in the context of percolation,
can be seen as random functions of the percolation configuration.
However, we view the graph as being fixed and satisfying the assumptions
stated in subsection \ref{sub:Assumptions-and-statement}.

We start with a general lemma on the relation between expected exit
time from a set and the expected number of visits to a fixed point
in the set.

Let $Z\subset\Gamma$ and let $\tau_{Z}$ be the first hitting time
of $Z$ for $X_{t}$. For $x\in\Gamma$ we set $G_{Z}(x)=E_{x}\left[\sum_{t=0}^{\tau_{Z}}\mathbf{1}_{\left\{ X_{t}=x\right\} }\right]$,
the expected number of visits to $x$ of a walk starting at $x$ before
$\tau_{Z}$.
\begin{lem}
There is a $k=k(d)>0$ such that for any $x\in\Gamma$ and $Z\subset\Gamma$
where $N(x)\cap Z=\emptyset$ \begin{equation}
E_{x}\left[\tau_{Z}\right]>kG_{Z}^{2}/\log G_{Z}.\label{eq:E_gt_G}\end{equation}
\end{lem}
\begin{proof}
Recall from \ref{sub:Assumptions-and-statement} that $X_{t}$ has
a positive staying probability at certain vertices. Since we apply
electrical network interpretation to estimate hitting probabilities,
we prove the above for $Y_{t}$, the usual discrete simple random
walk, with $0$ probability to stay at a vertex. This implies \eqref{eq:E_gt_G}
for $X_{t}$ as well, since the expected exit time cannot decrease,
and the Green functions for $X_{t}$ can only grow by a $2d$ factor
since the escape probabilites for $X_{t}$ and $Y_{t}$ are at most
a $2d$ factor apart.

Fix $x\in\Gamma$, set $T_{0}=0$, and define for each $i\in\mathbb{N}$
the r.v.'s \[
T_{i}=\inf\left\{ t>T_{i-1}:Y_{t}=x\right\} .\]
Let $i^{*}=\inf\left\{ i:T_{i}=\infty\right\} $. For $1\le i\le i^{*}$
let $\rho_{i}=T_{i}-T_{i-1}$. We show there are positive constants
$k_{1},k_{2}$ dependent only on $d$ such that 

\begin{equation}
P\left[\rho_{1}\ge k_{1}r^{2}/\log r\right]\ge\frac{k_{2}}{r}.\label{eq:excursion_len_low_bnd}\end{equation}
For some $r>1$, let $\partial B_{r}^{\Gamma}(x)=\left\{ v\in\Gamma:d_{\Gamma}(v,x)=r\right\} $.
By electrical network interpretation (see e.g. \citep{doyle1984random}),
the probability for a walk beginning at $x$ to hit $\partial B_{r}^{\Gamma}(x)$
before returning to $x$ is $(2d)^{-1}C_{eff}(r)$, where $C_{eff}(r)$
is the effective conductance from $x$ to $\partial B_{r}^{\Gamma}(x)$.
Since $\Gamma$ is infinite and connected, for any $r$ there is a
connected path of $r$ edges from $x$ to $\partial B_{r}^{\Gamma}(x)$.
By the monotonicity principle, $C_{eff}(r)$ is at least the conductance
on this path, which is $r^{-1}$. Thus the probability to hit some
$y\in\partial B_{r}^{\Gamma}(x)$ before returning to $x$ is at least
$\left(2dr\right)^{-1}$.

Next, let $y\in\partial B_{r}^{\Gamma}(x)$. By the Carne-Varopoulos
upper bound (see \citep{varopoulos1985long}), \[
P\left[Y_{t}=x|Y_{0}=y\right]\le4d^{1/2}\exp\left(-\frac{r^{2}}{2t}\right)\]
and thus, for some $k_{1}(d),k_{2}(d)>0$ and all $r>1$, the probability
that a walk starting at $y\in\partial B_{r}^{\Gamma}(x)$ does not
hit $x$ in the next $\lfloor k_{1}r^{2}/\log r\rfloor$ steps, by
union bound, is greater than $k_{2}$. Together with our lower bound
on the probability that we arrive at such a $y\in\partial B_{r}^{\Gamma}(x)$,
we get \eqref{eq:excursion_len_low_bnd}.

Next, let $g=\inf\left\{ i:T_{i}>\tau_{Z}\right\} $ be the number
of visits of $Y_{t}$ to $x$ before hitting $Z$, including $t=0$
. $g$ is a geometric random variable with mean $G=G_{Z}$. Let $\alpha=\frac{1}{2}\ln\left(4/3\right)$
and note that since there is a constant in \eqref{eq:E_gt_G} and
$\tau_{Z}\ge1$, we can assume $G>2$. Thus \begin{eqnarray*}
P\left[g\ge\alpha G\right] & \ge & \left(1-G^{-1}\right)^{\alpha G}\\
 & \ge & \left(1-G^{-1}\right)^{2\alpha\left(G-1\right)}\\
 & \ge & e^{-2\alpha}=3/4\end{eqnarray*}
We further assume $G>\frac{2}{\alpha}\vee\frac{16}{k_{2}}$ so that
$\alpha G-1>\alpha G/2$ and $G\frac{k_{2}}{16}>1$. Let $A$ be the
event that there is an $1\le i\le\left(\alpha G-1\right)\wedge i^{*}$
such that $\rho_{i}>k_{1}\left(\frac{k_{2}}{16}G\right)^{2}/\log\left(\frac{k_{2}}{16}G\right)$
Note that $i^{*}\le\alpha G-1$ implies $A$. Thus by \eqref{eq:excursion_len_low_bnd}
and the independence of consecutive excursions from $x$, \[
P\left[A^{c}\right]\le\left(1-\frac{16}{G}\right)^{\alpha G/2}\le e^{-8\alpha}\]
which is smaller than $1/4$.

Thus \[
P\left[\left\{ g\ge\alpha G\right\} ,A\right]\ge\frac{1}{2}.\]
This implies the lemma since $\tau_{Z}>\sum_{i=1}^{g-1}\rho_{i}$,
and for $k=k_{1}k_{2}^{2}/256$, $\sum_{i=1}^{g-1}\rho_{i}>kG^{2}/\log G$.
\end{proof}
Next, we state the fact that a Harnack inequality implies an oscillation
inequality. For a set of vertices $U$ and a function $u$, let 

\[
\mbox{osc}_{U}(u)=\max\nolimits _{U}(u)-\min\nolimits _{U}(u).\]

\begin{prop}
Let $x\in\Gamma$ and assume that for some $r>0$ we have that any
function $h$ that is non-negative and harmonic on $B_{2r}(x)$ satisfies
\eqref{eq:HK} on $B_{r}(x)$. Then for any $h$ that is harmonic
on $B_{2r}(x)$, we have\begin{equation}
\mbox{osc}_{B_{r}(x)}h\le\frac{c_{H}-1}{c_{H}+1}\mbox{osc}_{B_{2r}(x)}h.\label{eq:osc_ineq}\end{equation}
\end{prop}
\begin{proof}
We quote a proof from chapter 9 of \citep{telcs2006art}.

Set $v=h-\min_{B_{2r}(x)}h$. Since $v$ is non-negative and harmonic
on $B_{2r}(x)$, the Harnack inequality \eqref{eq:HK} is satisfied
here, so we have\[
\max_{B_{r}(x)}v\le c_{H}\min_{B_{r}(x)}v,\]
whence\[
\max_{B_{r}(x)}h-\min_{B_{2r}(x)}h\le c_{H}\left(\min_{B_{r}(x)}h-\min_{B_{2r}(x)}h\right),\]
and\[
\mbox{osc}_{B_{r}(x)}h\le\left(c_{H}-1\right)\left(\min_{B_{r}(x)}h-\min_{B_{2r}(x)}h\right).\]
Similarly, we have \[
\mbox{osc}_{B_{r}(x)}h\le\left(c_{H}-1\right)\left(\max_{B_{2r}(x)}h-\max_{B_{r}(x)}h\right).\]
Summing up these two inequalities, we obtain\[
\mbox{osc}_{B_{r}(x)}h\le\frac{1}{2}\left(c_{H}-1\right)\left(\mbox{osc}_{B_{2r}(x)}h-\mbox{osc}_{B_{r}(x)}h\right),\]
whence \eqref{eq:osc_ineq} follows.
\end{proof}
Iterating this on the Harnack assumption \ref{enu:HK} we get 
\begin{cor}
\label{cor:osc_iterate}For any $\alpha>0$ there is an $M(\alpha)$
such that for any $\eta>0$ and for all $R>R_{\eta}$, if $x\in B_{R}$,
$r\ge\eta R$ and $h$ is a harmonic function on $B_{M(\alpha)r}(x)$
then\textup{\[
\mbox{osc}_{B_{r}(x)}h\le\alpha\mbox{osc}_{B_{Mr}(x)}h.\]
}
\end{cor}
We use this to show regularity of the green function and for expected
exit time.
\begin{lem}
\label{lem:green_up_bnd}There is a $c_{G}(\epsilon)$ such that for
all large enough $R$, and any $x\in B_{(1-\epsilon)R}\setminus B_{\epsilon R}$\begin{equation}
G_{\tau_{R}}(\mathbf{0},x)<c_{G}R^{2-d}.\label{eq:green_up_bnd}\end{equation}
\end{lem}
\begin{proof}
Any two vertices in $u,v\in B_{(1-\epsilon)R}\setminus B_{\epsilon R}$
can be joined by a path of $n<C(d)/\epsilon$ overlapping balls $B_{\epsilon/4}(x_{1}),\ldots,B_{\epsilon/4}(x_{n})$
that are all subsets of $B_{(1-\epsilon/2)R}\setminus B_{\left(\epsilon/2\right)R}$
such that $x_{1}=u$ and $x_{n}=v$. Since $G_{\tau_{R}}(\mathbf{0},x)$
is positive and harmonic in $B_{(1-\epsilon/2)R}\setminus B_{\left(\epsilon/2\right)R}$,
Harnack assumption \ref{enu:HK} tells us that $c_{H}^{-C/\epsilon}<u/v<c_{H}^{C/\epsilon}$.
Next, note that $\sum_{x\in B_{R}}G_{\tau_{R}}(\mathbf{0},x)=E_{\mathbf{0}}\left[\tau_{R}\right]$.
Thus, letting $M(R)=\max\left\{ G_{\tau_{R}}(\mathbf{0},x):x\in B_{(1-\epsilon)R}\setminus B_{\epsilon R}\right\} $,
and using the exit time bound in assumption \ref{enu:exit_time_uni_bnd},
we get \[
\left|B_{R}\right|c_{H}^{-C/\epsilon}M\le\sum_{x\in B_{R}}G_{\tau_{R}}(\mathbf{0},x)<c_{E}R^{2}.\]
\end{proof}
\begin{lem}
\label{lem:Green_locally_lipschitz} For any $\beta>0$, there is
a $\delta>0$ such that for all large enough $R$, any $x,y\in B_{(1-\epsilon)R}\setminus B_{\epsilon R}$
satisfying $\left|x-y\right|<\delta R$ also satisfy\begin{equation}
\biggl|G_{\tau_{R}}(\mathbf{0},x)-G_{\tau_{R}}(\mathbf{0},y)\biggr|\le\beta R^{2-d}.\label{eq:Green_locally_Lip}\end{equation}
\end{lem}
\begin{proof}
$G_{\tau_{R}}(\mathbf{0},x)$ is positive and harmonic in $B_{(1-\epsilon/2)R}\setminus B_{\left(\epsilon/2\right)R}$
and bounded by $c_{G}R^{2-d}$. We use corollary \ref{cor:osc_iterate}
with  $\alpha=\beta c_{G}^{-1}$ and $\eta=\epsilon/\left(2M(\alpha)\right)$.
This gives us the lemma with $\delta=\eta$.\end{proof}
\begin{lem}
\label{lem:E_x-local-lip}For any $\beta>0$, there is a $\delta>0$
such that for all large enough $R$, any $x,y\in B_{(1-\epsilon)R}$
satisfying $\left|x-y\right|<\delta R$ also satisfy:
\end{lem}
\begin{equation}
\biggl|E_{x}\left[\tau_{R}\right]-E_{y}\left[\tau_{R}\right]\biggr|\le\beta R^{2}.\label{eq:E_x_local_lip}\end{equation}

\begin{proof}
Fix $x\in B_{(1-\epsilon)R}$ and for $\delta>0$ determined below,
fix some $y\in B_{\delta R}(x)$. Define the r.v. $\tau_{1}$ to be
the random time it takes a walk starting somewhere inside $B_{\epsilon R}(x)$
to exit $B_{\epsilon R}(x)$, and let $\tau_{2}=\tau_{R}-\tau_{1}$,
i.e. the additional time it takes the walk to exit $B_{R}$. Note
that $E_{x}\left[\tau_{2}\right]$ is harmonic in $B_{\epsilon R}(x)$,
and that from exit time assumption \ref{enu:exit_time_uni_bnd} for
all large $R$ it is bounded by $c_{E}R^{2}$. We use corollary \ref{cor:osc_iterate}
with $\alpha=\beta c_{E}^{-1}/2$ and $\eta=\epsilon/M(\alpha)$,
to get that for $\delta\le\eta$ \[
\biggl|E_{x}\left[\tau_{2}\right]-E_{y}\left[\tau_{2}\right]\biggr|\le\alpha c_{E}R^{2}=\beta R^{2}/2.\]
Take $\delta\le\left(\frac{\beta}{2c_{E}}\right)^{1/2}\wedge\eta$
so that, again by exit time bound, for all large enough $R$, $E_{x}\left[\tau_{1}\right],E_{y}\left[\tau_{1}\right]<\beta R^{2}/2$.
Applying the triangle inequality finishes the proof.
\end{proof}

\section{\label{sec:Domination-of-Green}Domination of Green function}

Let $\Omega=\Gamma^{\mathbb{N}\cup\left\{ 0\right\} }$ and let $P^{B}(\cdot)$
denote a probability measure on paths in $\Omega$ starting at $\mathbf{0}$
that make $X_{t}$ a {}``blind'' simple random walk as defined in
subsection \ref{sub:Assumptions-and-statement}. $P^{B}$ is pushed
forward to a measure on $\mathcal{C}_{T}^{d}$ by $w_{R}(t)$ as defined
in subsection \ref{sub:def-weak-conv-bm}. To contrast, we call $P(\cdot)$
the Wiener measure on curves corresponding to the Brownian motion
$\mathcal{B}(t)$ which is the weak limit of $w_{R}(t)$. Thus, for
fixed $T$, $\mathcal{C}_{T}^{d}$ is the probability space on which
$P^{B}\left(w_{R}(t)\right)$ converges to $P$ in distribution. Write
$E^{B}\left[\cdot\right]$ and $E\left[\cdot\right]$ for the corresponding
expectations.

\subsection{Integral Convergence of expected exit time}

Since we assume control of convergence to Brownian Motion only from
$\mathbf{0}$, we must describe the expected exit time from an arbitrary
point in the unit ball as a function of Brownian motion that starts
at $\mathbf{0}$. We do this by conditioning the Brownian motion to
hit a small box containing that point and measuring the additional
time needed to exit the unit ball.

We denote the first hitting time of a set $Z$ by $\tau_{Z}$. This
hitting times may refer to the Brownian motion $\mathcal{B}(t)$,
scaled linearly interpolated walks $w_{R}(t)$, or the discrete random
walk $X_{t}$. Another implicit part of this notation is the starting
point of the walk or curve. The correct interpretation should be evident
from context, and will be stated otherwise. Some notation used in
this section was introduced in subsections \ref{sub:Assumptions-and-statement}
and \ref{sub:def-weak-conv-bm}.

Fix $T>0$, $0<\theta<\epsilon$, $\mathbf{u}\in\mathbf{b}_{1-\epsilon}$
and let $\mathcal{A}=\mathcal{A}(T)=\left\{ w(t)\in\mathcal{C}_{2T}^{d}\ :\ \tau_{\mathbf{d}_{\theta}(\mathbf{u})}\le T\right\} $.
$\mathcal{A}\subset\mathcal{C}_{2T}^{d}$ is the event that the curve
hits a small box around $\mathbf{u}$ before time $T$. 

Henceforth, to avoid the complication where a vertex after scaling
is in the boundary of a box $\mathbf{d}_{\theta}(\mathbf{u})$, we
always take $\mathbf{u}$ to have rational coordinates, and the side
length $\theta$ to be irrational. This will suffice as our scale
parameter $R$ is a natural number. Secondly we take $T$ to be an
integer so that a curve $w_{R}(t)$ hits $\mathbf{d}_{\theta}(\mathbf{u})$
until $T$ if and only if it hits a vertex in $D_{R\theta}(R\mathbf{u})$
until $TR^{2}$.

Since we are interested in estimating the behavior of $X_{t}$ and
not just its interpolation, we define $\mathcal{A}^{*}=\mathcal{A}^{*}(R,T)=\left\{ \omega\in\Omega\ :\ \tau_{D_{R\theta}(R\mathbf{u})}(X_{t}(\omega))\le TR^{2}\right\} $.
$\mathcal{A}^{*}\subset\Omega$ is the event a vertex in $D_{R\theta}(R\mathbf{u})$
is visited until time $TR^{2}$. Thus for all $R\in\mathbb{N}$, and
for all $\omega\in\Omega$, $X_{t}(\omega)\in\mathcal{A}^{*}\iff w_{R}(\omega,t)\in\mathcal{A}$.%
{}

Let $\tau^{+}$ be the first exit time of the unit ball $\mathbf{b}_{1}$
after hitting $\mathbf{d}_{\theta}(\mathbf{u})$. Let $\mbox{\ensuremath{k:\mathcal{C}_{2T}^{d}\to\mathbb{R}}}$
be defined as \[
k(\omega)=\mathbf{1}_{\mathcal{A}}\bigl((\tau^{+}-\tau_{\mathbf{d}_{\theta}(\mathbf{u})})\wedge T\bigl).\]
Analogously, let $\tau_{R}^{+}$ be the first exit time of $B_{R}$
after hitting $D_{R\theta}(R\mathbf{u})$, and let \[
k^{*}(\omega)=R^{-2}\cdot\mathbf{1}_{\mathcal{A}^{*}}\bigl((\tau_{R}^{+}-\tau_{D_{R\theta}(R\mathbf{u})})\wedge TR^{2}\bigl).\]
$k(\omega)$ is bounded by $T$ and is discontinuous on $\partial\mathcal{A}$,
a set of Wiener measure zero. Therefore, by the Portmanteau theorem
and our assumption of weak convergence,\[
E^{B}\left[k\left(w_{R}(t)\right)\right]\to E\left[k\left(\mathcal{B}(t)\right)\right]\mbox{ as }R\to\infty.\]

By the strong Markov property for Brownian motion, we may average
over the hitting point.\[
P(\mathcal{A})^{-1}\cdot E\left[k\left(\mathcal{B}(t)\right)\right]=E\left[k\left(\mathcal{B}(t)\right)\ |\ \mathcal{A}\right]=E_{\mathbf{0}}\left[E_{\mathcal{B}(\tau_{\mathbf{d}_{\theta}(\mathbf{u})})}\left[\tau\wedge T\right]\right],\]
where $\tau=\tau_{\mathbf{b}_{1}^{c}}$, the first exit time from
the unit ball (we start measuring time at $\mathcal{B}(\tau_{\mathbf{d}_{\theta}(\mathbf{u})})$).
Note that $R^{2}k(w_{R}(\omega,t))$ measures the time that the unscaled
interpolated walk $w_{1}(\omega,t)$ takes to get from the boundary
of $\mathbf{d}_{R\theta}(R\mathbf{u})$ to $\mathbf{b}_{R}^{c}$,
but what interests us is the span between the first time that $X_{t}(\omega)$
takes a value in $D_{R\theta}(R\mathbf{u})$ to the first time it
takes a value in the complement of $B_{R}$. $R^{2}k^{*}\left(X_{t}(\omega)\right)$
measures this time. By the strong Markov property for random walks\[
P^{B}(\mathcal{A}^{*})^{-1}\cdot E^{B}\left[k^{*}\left(X_{t}\right)\right]=E^{B}\left[k^{*}\left(X_{t}\right)\ |\ \mathcal{A}^{*}\right]=R^{-2}E_{\mathbf{0}}^{B}\left[E_{Y(\tau_{D_{R\theta}(R\mathbf{u})})}^{B}\left[\tau_{R}\wedge TR^{2}\right]\right].\]
If the unscaled interpolated curve $w_{1}(\omega,t)$ crosses the
boundary of $\mathbf{d}_{R\theta}(R\mathbf{u})$, it will hit a vertex
in $D_{R\theta}(R\mathbf{u})$ in less than one time unit. The same
is true for exiting $\mathbf{b}_{R}.$ Thus for all $\omega$ in our
probability space $\left|R^{2}k\left(w_{R}(\omega,t)\right)-R^{2}k^{*}\left(X_{t}(\omega)\right)\right|<2$
and $\left|E^{B}\left[k(w_{R}(\omega,t))\right]-E^{B}\left[k^{*}(\omega)\right]\right|<\frac{2}{R^{2}}$.

By weak convergence, for any fixed $T,$ $P^{B}\left(w_{R}(\omega,t)\in\mathcal{A}(T)\right)\to P\left(\mathcal{A}(T)\right)$,
and since $\mathcal{A}^{*}\iff\mathcal{A}$,\begin{equation}
R^{-2}E_{\mathbf{0}}^{B}\left[E_{Y(\tau_{D_{R\theta}(R\mathbf{u})})}^{B}\left[\tau_{R}\wedge TR^{2}\right]\right]\to E_{\mathbf{0}}\left[E_{\omega(\tau_{\mathbf{d}_{\theta}(\mathbf{u})})}\left[\tau\wedge T\right]\right]\mbox{ as }R\to\infty.\label{eq:E_x-integral-conv}\end{equation}

In summary, we have some average on the boundary vertices of $D_{R\theta}(R\mathbf{u})$
of the function $R^{-2}E_{x}^{B}(\tau_{R}\wedge TR^{2})$, that is
close as we like to an average on $\mathbf{d}_{\theta}(\mathbf{u})$
of a Brownian motion's expected time to exit the unit ball.

\subsection{Integral convergence of Green function}

For a fixed $T>0$, $\theta>0$ and $\mathbf{u}\in\mathbf{b}_{1}$
let $h:\mathcal{C}_{T}^{d}\to\mathbb{R}$ be defined for $w(t)\in\mathcal{C}_{T}^{d}$
as follows:\[
h(w(t))=\int\limits _{0}^{T}\mathbf{1}_{\{w(t)\in\mathbf{d}_{\theta}(\mathbf{u}),t<\tau\}}dt.\]
$h(\omega)$ measures a curve's occupation time of $\mathbf{d}_{\theta}(\mathbf{u})$
before leaving $\mathbf{b}_{1}$ and until time $T$. Since $h(\mathcal{B}(\omega,t))$
is bounded by $T$ and is discontinuous on curves whose occupation
time of $\partial\mathbf{d}_{\theta}(\mathbf{u})$ before exit of
$\mathbf{b}_{1}$ is positive - a set of Wiener measure zero. Thus
the Portmanteau theorem gives:\[
E^{B}\left[h\left(w_{R}(t)\right)\right]\to E\left[h\left(\mathcal{B}(t)\right)\right]\mbox{ as }R\to\infty.\]
Note that $E\left[h\left(\mathcal{B}(t)\right)\right]=\int\limits _{\mathbf{d}_{\theta}(\mathbf{u})}g_{\tau\wedge T}(\mathbf{0},x)dx$
where $g_{\tau\wedge T}$ is the Green function of $\mathcal{B}(t)$
killed on leaving the unit ball or when time $T$ is reached.

Again, we are measuring the time that a linearly interpolated curve
spends in a set, while we would like to have control over the time
the random walk itself spends in the set. However, $E\left[h\left(\mathcal{B}(t)\right)\right]$
is a continuous function of $\theta,$ and for any $\delta>0$ the
time the discrete walk spends in $D_{R\theta}(R\mathbf{u})$ is eventually
sandwiched between the time the unscaled interpolated walk $w_{1}(\omega,t)$
spends in $\mathbf{d}_{R(\theta+\delta)}(R\mathbf{u})$ and $\mathbf{d}_{R(\theta-\delta)}(R\mathbf{u})$.
Thus:

\begin{equation}
R^{-2}\sum_{x\in D_{R\theta}(R\mathbf{u})}G_{\tau_{R}\wedge TR^{2}}^{B}(\mathbf{0},x)\to\int\limits _{\mathbf{d}_{\theta}(\mathbf{u})}g_{\tau\wedge T}(\mathbf{0},x)dx\mbox{ as }R\to\infty,\label{eq:Green_integral_convergence}\end{equation}
where $G_{\tau_{R}\wedge TR^{2}}^{B}$ is the Green function of the
walk from $\mathbf{0}$, killed on leaving $B_{R}$ or when time $TR^{2}$
is reached.

\subsection{Pointwise domination of Green function}
\begin{lem}
\label{lem:pointwise_domination}$\forall\epsilon>0\ \exists\hat{R}\mbox{ and }\eta>0\mbox{ s.t. }\forall R>\hat{R}$
the following holds: \[
x\in B_{(1-\epsilon)R}\setminus B_{\epsilon R}\implies|B_{R}|G_{\tau_{R}}(\mathbf{0},x)>(1+\eta)E_{x}(\tau_{R}).\]

\end{lem}
This is the main result needed for IDLA lower bound, and will proved
be in this section.

First, it is known (see, e.g., \citep{lawler1992internal} p.2121)
for Brownian motion starting at zero, and killed on exiting the unit
ball, that the Green function $g_{\tau}(0,x)$ and expected hitting
time $E_{x}(\tau)$ are continuous functions with the property that
$|\mathbf{b}_{1}|g_{\tau}(\mathbf{0},x)-E_{x}(\tau)$ descends strictly
monotonically to zero as $x$ goes from $0$ to $1$. This is true
for any Brownian motion, in particular, $\mathcal{B}(t)$, the weak
limit of $X_{t}$ on the graph starting at $\mathbf{0}$. Thus for
any $\epsilon>0$ the minimum of the difference between the two when
$||x||\in[\epsilon,1-\epsilon]$ is bounded away from zero. The Lebesgue
monotone convergence theorem implies that $g_{\tau\wedge T}(\mathbf{0},x)\nearrow g_{\tau}(\mathbf{0},x)$
and $E_{x}(\tau\wedge T)\nearrow E_{x}(\tau)$ as $T\to\infty$. Since
all functions involved are continuous and converge monotonically on
a compact set, by Dini's theorem, the convergence is uniform. Thus
we have the following uniform bounding of the difference away from
zero: \[
\exists\gamma(\epsilon)>0,\ T\mbox{ s.t. }\forall\mathbf{u},||\mathbf{u}||\in[\epsilon/2,1-\epsilon/2],\ E_{\mathbf{u}}(\tau\wedge T)+\gamma<|\mathbf{b}_{1}|g_{\tau\wedge T}(\mathbf{0},\mathbf{u}).\]

Since $E_{\mathbf{u}}(\tau\wedge T)$, $g_{\tau\wedge T}(0,\mathbf{u})$
converge uniformly with $T$, they are uniformly equicontinuous in
the variable $\mathbf{u}$ in the closed interval $[\epsilon/2,1-\epsilon/2]$.
%
{}We may then choose a $\theta>0$ such that for all large enough $T$,
any average of $g_{\tau\wedge T}(\mathbf{0},\cdot)$ in a box of side
$\theta$ with center $\mathbf{u}$, is close to $g_{\tau\wedge T}(0,\mathbf{u})$
for any $\mathbf{u\in[\epsilon,1-\epsilon]}$. We have the analogous
claim for $E_{\mathbf{u}}(\tau\wedge T)$. Thus:
\begin{lem}
\label{lem:Domination-of-cont-green-in-small-ball}For any positive
$\epsilon$, there is $\gamma(\epsilon)>0$ such that for all large
enough $T$ and all small enough $\theta$\[
\ \forall\mathbf{u},||\mathbf{u}||\in[\epsilon,1-\epsilon],\ \int\limits _{\partial\mathbf{d}_{\theta}(\mathbf{u})}E_{\mathbf{v}}(\tau\wedge T)d\mu(\mathbf{v})+3\gamma<|\mathbf{b}_{1}|\theta^{-d}\int\limits _{\mathbf{d}_{\theta}(\mathbf{u})}g_{\tau\wedge T}(\mathbf{0},\mathbf{v})d\mathbf{v}.\]

\end{lem}
In the above, $d\mu$ is an arbitrary probability measure (total mass
one) on the boundary of $\mathbf{d}_{\theta}(\mathbf{u})$, while
the integral on the right is by $d$-dimensional Lebesgue measure.

We apply lemmas \ref{lem:Green_locally_lipschitz} and \ref{lem:E_x-local-lip}
to get a $\delta>0$ for which \eqref{eq:E_x_local_lip} and \eqref{eq:Green_locally_Lip}
hold with $\beta=\hat{\gamma}$ for all large enough $R$. $\hat{\gamma}>0$
is some multiple of $\gamma$ from lemma \ref{lem:Domination-of-cont-green-in-small-ball}
that we determine later. Set $\theta$ to be small enough so that
$\mathbf{d_{\theta}}$ is covered by a ball of radius $\delta$. Increase
$T$ further if necessary so that \eqref{eq:E_x-uni-conv-with-T}
holds with $\beta=\frac{\gamma}{4}$. Fix $T$ and $\theta$ for the
remainder of the proof.

We cover $\mathbf{b}_{(1-\epsilon)}\setminus\mathbf{b}_{\epsilon}$
by a finite number of open $\theta$-boxes, and so to prove lemma
\ref{lem:pointwise_domination}, it suffices to prove the implication
in the restricted setting of $x\in D_{R\theta}(R\mathbf{u})$ where
$\mathbf{u}$ is the center of an arbitrary box in our $\theta$-net.

Now, we show the Green function at every point in this box is close
to the continuous one.

\subsubsection{Pointwise Green function estimate}

Let $G_{R}^{\Sigma}=\sum\limits _{x\in D_{R\theta}(R\mathbf{u})}G_{\tau_{R}\wedge TR^{2}}(\mathbf{0},x)$
and let $g^{\int}=\int\limits _{\mathbf{d}_{\theta}(\mathbf{u})}g_{\tau\wedge T}(\mathbf{0},x)dx$.
By \eqref{eq:Green_integral_convergence}, we have for large enough
$R$:\[
\biggl|R^{-2}G_{R}^{\Sigma}-g^{\int}\biggr|<\hat{\gamma}\theta^{d}|\mathbf{b}_{1}|.\]
Let $\alpha_{R}=\frac{|B_{R}|}{|D_{R\theta}(R\mathbf{u})|}$. By density
assumption \eqref{enu:Convergence-to-positive}, for all $R$ large
enough $\biggl|\alpha_{R}-\theta^{-d}|\mathbf{b}_{1}|\biggr|<\hat{\gamma}$.
Using the triangle inequality on \[
\biggl|R^{-2}G_{R}^{\Sigma}(\alpha_{R}+\theta^{-d}|\mathbf{b}_{1}|-\alpha_{R})-\theta^{-d}|\mathbf{b}_{1}|g^{\int}\biggr|<\hat{\gamma},\]
we get\[
\biggl|R^{-2}G_{R}^{\Sigma}\alpha_{R}-\theta^{-d}|\mathbf{b}_{1}|g^{\int}\biggr|<\hat{\gamma}+R^{-2}G_{R}^{\Sigma}\Bigl|\alpha_{R}-\theta^{-d}|\mathbf{b}_{1}|\Bigr|<\hat{\gamma}(1+c_{5}\theta^{d}),\]
where for the right inequality, we used \eqref{eq:green_up_bnd} to
bound $G_{\tau_{R}\wedge TR^{2}}(\mathbf{0},x)$ and $(R\theta)^{d}$
as a bound on the number of vertices.

We set $\theta$ such that \eqref{eq:Green_locally_Lip} bounds the
difference between the maximum and minimum of $G_{\tau_{R}}(\mathbf{0},x)$
in $D_{R\theta}(R\mathbf{u})$ by $\hat{\gamma}R^{2-d}$ for all large
enough $R$. Thus, since $G_{\tau_{R}}(\mathbf{0},x)>G_{\tau_{R}\wedge TR^{2}}(\mathbf{0},x)$
we have for any $x\in D_{R\theta}(R\mathbf{u})$ \[
\ \left(R^{-2}G_{\tau_{R}}(\mathbf{0},x)|B_{R}|-\theta^{-d}|\mathbf{b}_{1}|g^{\int}\right)>-\left(\hat{\gamma}(1+c_{5}\theta^{d})+R^{-2}|B_{R}|\hat{\gamma}R^{2-d}\right).\]
Recall $\gamma$ from lemma \ref{lem:Domination-of-cont-green-in-small-ball}.
We now determine $\hat{\gamma}$ so that for all large enough $R$,
for any $x\in D_{R\theta}(R\mathbf{u})$ \begin{equation}
\left(R^{-2}G_{\tau_{R}}(\mathbf{0},x)|B_{R}|-\theta^{-d}|\mathbf{b}_{1}|\int\limits _{\mathbf{d}_{\theta}(\mathbf{u})}g_{\tau\wedge T}(\mathbf{0},\mathbf{v})d\mathbf{v}\right)>-\gamma.\label{eq:pt-discrete-green-close-to-continuous}\end{equation}
Next we show that the expected hitting time from any point is close
to the continuous one.

\subsubsection{Pointwise expected hitting time estimate}

Recall we chose $T$ such that for any $x\in B_{R}$ \[
\left|E_{x}\left[\tau_{R}\right]-E_{x}\left[\tau_{R}\wedge TR^{2}\right]\right|<\frac{\gamma}{4}R^{2}.\]
Combining this with \eqref{eq:E_x-integral-conv} we have for large
enough $R$ \[
\Biggl|R^{-2}\sum_{x\in D_{R\theta}(R\mathbf{u})}\lambda_{x}^{(R)}E_{x}(\tau_{R})-\int\limits _{\partial\mathbf{d}_{\theta}(\mathbf{u})}E_{\mathbf{v}}(\tau\wedge T)d\mu(\mathbf{v})\Biggr|<\frac{\gamma}{2}.\]
 where for each $R$, $\lambda_{x}^{(R)}$ are non-negative and sum
to one and $d\mu$ is some probability measure on $\partial\mathbf{d}_{\theta}(\mathbf{u})$. 

We set $\theta$ such that \eqref{eq:E_x_local_lip} holds with $\beta=\hat{\gamma}<\gamma/2$,
so for large enough $R$ we have for any $x\in D_{R\theta}(R\mathbf{u})$
\[
\Biggl|R^{-2}E_{x}(\tau_{R})-\int\limits _{\partial\mathbf{d}_{\theta}(\mathbf{u})}E_{\mathbf{v}}(\tau\wedge T)d\mu(\mathbf{v})\Biggr|<\gamma.\]
Putting the above together with \eqref{eq:pt-discrete-green-close-to-continuous}
and \ref{lem:Domination-of-cont-green-in-small-ball}, for all large
enough $R$, we get that for any $x\in D_{R\theta}(R\mathbf{u})$\[
G_{\tau_{R}}(\mathbf{0},x)|B_{R}|-E_{x}(\tau_{R})>\gamma R^{2}.\]
The above along with our upper bound on $E_{x}(\tau_{R})$ from assumption
\ref{enu:exit_time_uni_bnd}, implies lemma \ref{lem:pointwise_domination}.

\section{\label{sec:Lower-Bound}Lower Bound}

We begin by formally defining the IDLA process.

Let $\left(X_{t}^{n}\right)_{t\ge0}^{n\in\mathbb{N}}$ be a a sequence
of independent random walks starting at $\mathbf{0}$. The aggregate
begins empty, i.e. $\mathcal{I}(0)=\emptyset$, and $\mathcal{I}(n)=\mathcal{I}(n-1)\cup X_{t'}^{n}$
where $t'=\min_{t}\left\{ X_{t}^{n}\notin\mathcal{I}(n-1)\right\} $.
Thus we have an aggregate growing by one vertex in each time step.

As in \citep{lawler1992internal}, we fix $z\in B_{(1-\epsilon)R}\setminus B_{\epsilon R}$,
and look at the first $\left|B_{R}\right|$ walks. Let $A=A(z,R)$
be the event $z\in\mathcal{I}(\left|B_{R}\right|)$. We show the probability
this does not happen decreases exponentially with $R$. 

Let $M=M(z,R)$ be the number of walks out of the first $\left|B_{R}\right|$
that hit $z$ before exiting $B_{R}$. Let $L=L(z,R)$ be the number
of walks out of the first $\left|B_{R}\right|$ that hit $z$ before
exiting $B_{R}$, but after leaving the aggregate. Then for any $a$,
\[
P(A^{c})<P(M=L)<P(M\le a)+P(L\ge a).\]
We choose $a$ later to minimize the terms. In order to bound the
above expression, we calculate the average of $M$ and $L$, 

\[
E\left[M\right]=\left|B_{R}\right|P_{\mathbf{0}}(\tau_{z}<\tau_{R}).\]
$E\left[L\right]$ is hard to determine, but each walk that contributes
to $L$, can be tied to the unique point at which it exits the aggregate.
Thus, by the Markov property, if we start a random walk from each
vertex in $B_{R}$ and let $\hat{L}$ be those walks that hit $z$
before exiting $B_{R}$, $P(L>n)\le P(\hat{L}>n)$, and so it suffices
to bound $P(\hat{L}>a)$. $E[\hat{L}]$ is a sum of independent indicators:\[
E[\hat{L}]=\sum_{x\in B_{R}}P_{x}(\tau_{z}<\tau_{R}).\]

Since both $M$ and $\hat{L}$ are sums of independent variables,
we expect them to be close to their mean. Our aim now becomes showing
that for some $\delta>0$ and all large enough $R$,\begin{equation}
E\left[M\right]>(1+2\delta)E[\hat{L}].\label{eq:E_M-gt-E_L}\end{equation}
By standard Markov chain theory,\[
P_{x}(\tau_{z}<\tau_{R})=\frac{G_{\tau_{R}}(x,z)}{G_{\tau_{R}}(z,z)}.\]
Using this and symmetry of the Green function, it is enough to show
for all large enough $R$\[
\left|B_{R}\right|G_{\tau_{R}}(\mathbf{0},z)>(1+2\delta)\sum_{x\in B_{R}}G_{\tau_{R}}(x,z)=(1+2\delta)E_{z}\left[\tau_{R}\right],\]
which is lemma \ref{lem:pointwise_domination}.

Choosing $a=(1+\delta)E[\hat{L}]$ we write

\begin{eqnarray}
P(\hat{L}>(1+\delta)E[\hat{L}]) & < & P(\Bigl|\hat{L}-E[\hat{L}]\Bigr|>\delta E[\hat{L}]^{1/2}\sigma_{\hat{L}})\\
 & < & 2\exp\left(-E[\hat{L}]\delta^{2}/4\right).\label{eq:L_hat_bound}\end{eqnarray}
In the first line we use that $\hat{L}$ is a sum of independent indicators,
and the variance of such a sum is smaller than the mean. The second
line is an application of Chernoff's inequality. Similarly, using
\eqref{eq:E_M-gt-E_L} 

\begin{eqnarray}
P(M<(1+\delta)E[\hat{L}]) & < & P(M<\frac{1+\delta}{1+2\delta}E\left[M\right])\\
 & < & P(\left|M-E\left[M\right]\right|>\frac{\delta}{2}E\left[M\right])\\
 & < & 2\exp\left(-E\left[M\right]\delta^{2}/16\right).\label{eq:M-bound}\end{eqnarray}
To lower bound $E[\hat{L}]=E_{z}[\tau_{R}]\left(G_{\tau_{R}}(z,z)\right)^{-1}$
we use \eqref{eq:E_gt_G} to write that for some $M$,\[
G_{\tau_{R}}(z,z)\le\left(ME_{z}[\tau_{R}]\log E_{z}[\tau_{R}]\right)^{1/2}\]
and since $E_{z}[\tau_{R}]\ge R$ we have \[
E\left[\hat{L}(z,R)\right]>\frac{cR^{1/2}}{\log R}.\]
Together with \eqref{eq:L_hat_bound} and \eqref{eq:M-bound} and
summing over all $z\in B_{(1-\epsilon)R}\setminus B_{\epsilon R}$,
we get that the probability one of the vertices in $B_{(1-\epsilon)R}\setminus B_{\epsilon R}$
is not in $\mathcal{I}(\left|B_{R}\right|)$ is bounded by $CR^{d}\exp(-cR^{1/2}/\log R)$.
Since the expression is summable by $R$, by Borel Cantelli, this
happens only a finite number of times with probability one. So if
$R'$ is the largest radius for which some vertex in $B_{R'(1-\epsilon)}\setminus B_{R'\epsilon}$
is not covered after $\left|B_{R'}\right|$ steps, then we possibly
have a finite sized hole in the aggregate which will almost surely
fill up after another finite number of steps.

Thus we have proved the main theorem \ref{thm:IDLA-lower-bound}.

\emph{Acknowledgement. }Thanks to Itai Benjamini for suggesting the
problem, to Gady Kozma for many helpful discussions, and to Greg Lawler
for suggesting an elegant approach.

\bibliographystyle{amsalpha}
\bibliography{IDLASCC}

\end{document}